\documentclass[draft,12pt]{article}
\usepackage{amsmath,amsthm}
\newtheorem{thm}{Theorem}

\newtheorem{cor}{Corollary}
\newtheorem{lem}{Lemma}

\theoremstyle{remark}

\numberwithin{equation}{section}
\def\Pr{{\rm Pr}}

\def\epsilon{\varepsilon}

\begin{document}

\title{A new random mapping model
 \\}

\author{Jennie C. Hansen\footnote{Actuarial Mathematics and Statistics,
Heriot--Watt University,Edinburgh EH14 4AS, UK.
E-mail address:
J.Hansen@ma.hw.ac.uk}  ~and Jerzy Jaworski\footnote{Faculty of Mathematics and Computer Science, Adam Mickiewicz University, Umultowska 87, 61-614 Pozna\'n, Poland.
E-mail address: jaworski@amu.edu.pl
}\footnote{J.Jaworski
acknowledges the generous support by the
Marie Curie Intra-European Fellowship No. 501863 (RANDIGRAPH) within the 6th
 European Community Framework Programme.} }

\maketitle

\begin{abstract}
In this paper we introduce a new random mapping model, $T_n^{\hat D}$, which
maps the set $\{1,2,...,n\}$ into itself.
The random mapping $T_n^{\hat D}$ is constructed using a collection of exchangeable random variables $\hat{D}_1, ....,\hat{D}_n$ which satisfy $\sum_{i=1}^n\hat{D}_i=n$. In the random
digraph,
$G_n^{\hat D}$, which represents the mapping $T_n^{\hat D}$,
the in-degree sequence for the vertices is given by the variables $\hat{D}_1, \hat{D}_2, ..., \hat{D}_n$, and, in some sense,
$G_n^{\hat D}$
can be viewed as an analogue of the general independent degree
models  from random graph theory. We show that the distribution of the number of cyclic points, the number of components,
and the size of a typical component  can be
expressed in terms of  expectations of various functions of $\hat{D}_1, \hat{D}_2, ..., \hat{D}_n$. We also consider two special examples of $T_n^{\hat D}$
 which correspond to random mappings with preferential
and anti-preferential attachment, respectively, and determine, for these examples, exact and asymptotic distributions for the statistics mentioned above.
Results for the distribution of the number of successors and predecessors of a typical vertex in $G_n^{\hat D}$ in terms of expectations of various functions of $\hat{D}_1, \hat{D}_2, ..., \hat{D}_n$ are obtained in a companion paper
\cite{HJ06a}. \end{abstract}

\section{Introduction }
\label {s:intro} The study of random mapping models was initiated
independently by several authors
in the 1950s 
(see  \cite{Austin59, Folkert55, Ford57, Harris60, Katz55, Rubin54}
) and the properties of these models have received much attention in
the literature. In particular, these models have been useful as
models for epidemic processes,
 and have natural applications in cryptology (see, for example, \cite{Ball97, Berg81, Berg83, Burtin80,
 Flajolet90, Gertsbakh77, Jaworski98, Jaworski99, Mutafchiev81, Mutafchiev82, Pittel83, Qu90, VW94}).  To date, the most widely studied models have
 been special cases of a general model
denoted by $T_{{\bf p}(n)}$, which can be defined as follows: Let $[n]$ denote the set of integers
 $\{1,2,...,n\}$ and let ${\cal M}_n$ denote the set of all mappings from $[n]$ into $[n]$.
For each $n\geq 1$, let ${\bf p}(n)=\{p_{ij}(n) :1\leq i,j\leq n\}$ be an array
 such that $p_{ij}(n)\geq 0$ for $1\leq i,j\leq n$ and $\sum_{j=1}^n p_{ij} (n)=1$
 for every $1\leq i\leq n$, and let $X_1^n, X_2^n, ..., X_n^n$ be independent random variables
 such that $\Pr\{X_i^n=j\}=p_{ij}(n)$ for all $1\leq i,j\leq n$. Then the random mapping
 $T_{{\bf p}(n)}:[n]\to[n]$ is defined (in terms of the variables $X_1^n, X_2^n, ..., X_n^n$)
 by
\begin{equation}
\label{defn}
T_{{\bf p}(n)}(i)=j \phantom{X}{\rm iff}\phantom{X}
X_i^n=j
\end{equation}
for all $1\leq i,j\leq n$.
It follows from (\ref{defn}) that the distribution of $T_{{\bf
p}(n)}$ is given by
 \begin{equation}
 \label{dist}
\Pr\bigl\{ T_{{\bf p}(n)}=f \bigr\}=\prod_{i=1}^n p_{if(i)}(n)
\end{equation}
for each $f\in{\cal M}_n$. Any mapping $f\in{\cal M}_n$ can be represented
as a directed graph $G(f)$ on a set of vertices labelled $1,2,..., n$, such that
there is a directed edge from vertex $i$ to vertex $j$ in $G(f)$ if and only if
$f(i)=j$. So $G_{{\bf p}(n)}\equiv G(T_{{\bf p}(n)})$ is a random directed
graph on a set of vertices labelled $1,2,..., n$ which represents the action of the
random mapping
$T_{{\bf p}(n)}$ on $[n]$.
 We note that since each
vertex in $G_{{\bf p}(n)}$ has out-degree 1, the components of $G_{{\bf p}(n)}$
consist of directed cycles with directed trees attached. Also, it follows from the definition of
$T_{{\bf p}(n)}$ that the
 variables $X_1^n, X_2^n, \ldots X_n^n$ can be interpreted as
 the independent `choices' of the vertices $1,2, \ldots ,n$
in the random digraph $G_{{\bf p}(n)}$ (see, in addition,  Mutafchiev \cite{Mutafchiev84} and Jaworski \cite{Jaworski90}).

\par
The  example of $T_{{\bf p}(n)}$ which is best understood is
the uniform random mapping, $T_n\equiv T_{{\bf p}(n)}$, where $p_{ij}(n)={1\over n}$ for all
$1\leq i,j\leq n$ . Much is known (see for example the monograph by Kolchin
\cite{Kolchin86}) about the component structure of the
random digraph $G_n\equiv G(T_n)$ which represents $T_n$.
Aldous \cite{Aldous85} has shown that the
joint distribution of the normalized order statistics for the component
sizes in $G_n$ converges to the {\it Poisson-Dirichlet} (1/2) distribution
on the simplex $\nabla=\{\{x_i\}:\sum x_i\leq 1, x_i\geq x_{i+1}\geq 0$
for every $i\geq 1\}$.
Also, if $M_k$ denotes the number of components of size
$k$ in $G_n$ then the joint distribution
of $(M_1, M_2, \ldots , M_b)$ is close, in the sense of total variation,
to the joint distribution of a sequence
of independent Poisson random variables when $b=o(n/\log n)$
(see Arratia et.al. \cite{Arratia92}, \cite{Arratia95}) and from this result one obtains a
functional central limit theorem for the component sizes (see also \cite{Hansen89}).
The asymptotic distributions
of variables such as the number of predecessors and the
number of successors of a vertex in $G_n$ are also known (see
\cite{ Berg81, Berg83, Burtin80, Jaworski98, Jaworski99,
Mutafchiev81, Mutafchiev82, Pittel83}).
In another direction, Berg, Jaworski, and Mutafchiev (see \cite{BergJa89, Jaworski90, JM94,  Jaworski98,  Jaworski99} )
have investigated the structure of
$G_{{\bf p}(n)}$ when ${{\bf p}(n)}$ is given by $p_{ii}(n)=q$ for some
$0\leq q\leq 1$ and all $1\leq i\leq n$, and $p_{ij}(n)={1-q\over n-1}$
for all $1\leq i,j\leq n$ such that $i\not=j$.
Finally,   Aldous, Miermont, and Pitman (see \cite{AMP04}
and \cite{AMP05}) have recently  investigated the asymptotic structure of $G_{{\bf p}(n)}$,
where ${\bf p}(n)$ is given by $p_{ij}(n)=p_j(n)>0$ for all $1\leq i,j\leq n$, by using
an ingenious coding of the mapping $T_{{\bf p}(n)}$ as a stochastic process on the interval
$[0,1]$. Their results are closely related to earlier work on the relationship between
random mappings and random forests (see Pitman \cite{Pitman01} and references therein).

The common feature in all the models discussed above is that each vertex in $G_{{\bf p}(n)}$
`chooses' the vertex that it is mapped to independently of the `choices' made by all other vertices.
In this paper we introduce a new random mapping model in which the vertex
`choices' are not necessarily independent. The definition of the model is motivated, in part,
by  developments in the general theory of random graphs.
In recent years models for random graphs with
a specified degree sequence have received much attention as models for
complex networks such as the internet.
Loosely speaking, such a random graph on $n$ labelled
vertices can be constructed by starting with a collection of i.i.d., non-negative,
integer-valued random variables $D_1, D_2, ..., D_n$ and adding edges, at random, to the
graph until each vertex $i$ has degree $D_i$ in the constructed random graph.
Of interest in such models is the relationship between the component structure of the
random graph and  the distribution
of the variables $D_1, D_2,..., D_n$.  In another direction random graph models with `preferential attachment'
have been constructed  in order to model the {\it evolving} structure of
complex networks. In such models edges are added sequentially to the graph
and new edges are more likely to be attached to vertices that already have relatively high
degree in the evolving graph. The literature on these new developments in the theory
of random graphs is extensive, but a good bibliography is provided by Bonato's
survey paper \cite{Bonato}.

 \par In this paper we consider a
random mapping analogue for both the independent degree models and
the preferential attachment models
 described above. We show that in the case of random mappings these analogues are equivalent and we develop a  calculus for determining the distributions of various
important random mapping statistics which is based on the underlying distribution of the in-degrees
of the vertices in the directed graph which represents  the random mapping.
The paper is organized as follows. In Section 2 we define carefully define
the new random mapping model. In Section 3 we derive the calculus for this new model. In Section 4
we define both a random mapping model with preferential attachment and with
anti-preferential attachment. We show that both of these models are equivalent,
under certain distribution assumptions, to the special examples of the general model defined in Section 2.
For these special examples we investigate the distribution of the
number of cyclic points, the distribution of the number of components,
the probability of connectedness, and the distribution
of the size of a typical component.

\section{The model }
\label {s:model}
In order to define our new random mapping model, we adopt the following notation.
For $n\geq 1$,   ${\cal M}_n$
denotes the set of all mappings $f:[n]\to [n]$,  where
$[n]\equiv\{1,2,...,n\}$, and
$G(f)$ denotes, as described in the Introduction, the directed graph on $n$ labelled vertices
which represents the mapping $f\in{\cal M}_n$.
In addition, for $1\leq i\leq n$,  $d_i(f)$ denotes the in-degree of vertex $i$
in the digraph $G(f)$, and we let $\vec{d}(f)\equiv (d_1(f), ..., d_n(f))$.
\par
Now suppose that $\hat{D}_1, ....,\hat{D}_n$ is a collection of
exchangeable random va\-riables such that $\sum_{i=1}^n\hat{D}_i=n$, then
we construct a probability measure
${\bf P}_n^{\hat{D}}$ on ${\cal M}_n$ as follows.  For $f\in {\cal M}_n$,
define
\begin{equation}
\label{meas}
{\bf P}_n^{\hat D}\{f\}=\frac{\prod_{i=1}^n (d_i(f))!}{n!}
\Pr\big\{\hat{D}_i=d_i(f), 1\leq i\leq n  \big\}.
\end{equation}
It is clear from the definition of ${\bf P}_n^{\hat{D}}$ that  for any $f,g\in{\cal M}_n$ such that
$\vec{d}(f)=\vec{d}(g)$, we have ${\bf P}_n^{\hat{D}}\{f\}={\bf P}_n^{\hat{D}}\{g\}$.
Given the probability measure ${\bf P}_n^{\hat{D}}$, we can define the random mapping $T_n^{\hat{D}}$ which  takes
 values in ${\cal M}_n$ and has distribution
given by
\begin{equation}
\label{dist}
\Pr\big\{T_n^{\hat{D}}=f\big\}={\bf P}_n^{\hat{D}}\{f\}
\end{equation}
for every $f\in{\cal M}_n$, and we let $G_n^{\hat{D}}\equiv G(T_n^{\hat{D}})$ denote the
random digraph on $n$ labelled vertices which represents $T_n^{\hat{D}}$.
\par
An important class of examples can be constructed as follows.
Let $D_1, $ 
$D_2, \ldots, D_n$ be i.i.d.
non-negative integer-valued random variables. It does not make sense to construct directly
a random mapping digraph with in-degrees given by $D_1, D_2, ..., D_n$
since the sum of the vertex in-degrees in a random mapping digraph on $n$ vertices
always equals $n$. So, instead, we let
 $\hat{D}_1, \hat{D}_2,...,\hat{D}_n$ be a sequence of random variables such that
 their joint distribution is given by
$$
\Pr\big\{\hat{D}_i=d_i, 1\leq i\leq n\big\}=\Pr\Big\{D_i=d_i, 1\leq i\leq n \,\Big| \,\sum_{i=1}^n
D_i=n\Big\}
$$
and we use $\hat{D}_1, \hat{D}_2, ..., \hat{D}_n$ to construct $T_n^{\hat D}$ and $G_n^{\hat D}$.
This gives us the natural analogue of the general  i.i.d. degree model discussed in the
Introduction.
We also remark that it is easy to check that if $D_1, D_2, ..., D_n$
are i.i.d. $ Poisson(1)$ variables, then
the corresponding random mapping $T_n^{\hat D}$ is just the usual uniform
random mapping.
We will
see in Section 4 below that there are interesting interpretations $T_n^{\hat D}$
 in the cases where the underlying i.i.d. variables
 $D_1, D_2, ..., D_n$ have {\bf (i)}
 a generalised negative binomial distribution, and
{\bf (ii)}
a $ Bin(m, 1/m)$ distribution
for some $m\geq 2$.

\section{Results}
\label {s:results}
 In this section we develop a calculus, in terms of  the variables $\hat{D}_1, \hat{D}_2, ...,
\hat{D}_n$, for determining the distributions
of various random variables associated with the structure of
$G_n^{\hat{D}}$ .
The first variable we consider is the number of cyclic
vertices in the random digraph $G_n^{\hat{D}}$.
\par A vertex $i\in [n]$ is a cyclic vertex
for the mapping $f\in{\cal M}_n$ (and for the corresponding digraph $G(f)$
if there is some $k\geq 1$ such that $f^{(k)}(i)=i$, where $f^{(k)}$
is the $k^{th}$ iterate of the function $f$. We define $X_n(f)$ to be
the number of cyclic vertices of $f\in{\cal M}_n$ and we let $X_n^{\hat{D}}\equiv X_n(T_n^{\hat{D}})$
denote the number of cyclic vertices in $G_n^{\hat{D}}$. Then
we have

\begin{thm}
\label{Theorem1}
For $1\leq k\leq n-1$
$$
\Pr\big\{X_n^{\hat{D}}=k \big\}=\frac{k}{n-k}E\big((\hat{D}_1-1)\hat{D}_1\hat{D}_2\cdot\cdot\cdot \hat{D}_k\big)
$$
and
$$
\Pr\big\{X_n^{\hat{D}}=n \big\}= \Pr\big\{\hat{D}_i=1, 1\leq i\leq n \big\}.
$$
\end{thm}
\begin{proof} {\it (of Theorem 1)\,} We begin by considering the case $1\leq k\leq n-1$. For $f\in {\cal M}_n$,
let ${\cal L}_n(f)$ denote the set of cyclic vertices for the mapping $f$, and let
${\cal L}_n^{\hat{D}}\equiv {\cal L}_n(T_n^{\hat{D}})$. Then we have
$$
\Pr\{X_n^{\hat{D}}=k\}=\sum_{L\subseteq [n] s.t. |L|=k}\Pr\{{\cal L}_n^{\hat{D}}=L\}.
$$
We fix $L=\{1,2,...,k\}=[k]$ and observe that
$$
\Pr\big\{{\cal L}_n^{\hat{D}}=[k]\big\}=
$$
$$
\sum_{{\{d_i\}}\atop  {s.t.\sum d_i=n}}\Pr\Big\{{\cal L}_n^{\hat{D}}=[k] \, \Big|\,{\hat{D}}_i=d_i, 1\leq i\leq n
\Big\}
\Pr\big\{ {\hat{D}}_i=d_i, 1\leq i\leq n\big\}.
$$
Now fix  $\vec{d}=(d_1, d_2,..., d_n)$ such that $\sum_{i=1}^n d_i=n$
and let
$$
{\cal M}_n(\vec{d})=
\{f\in{\cal M}_n : d_i(f)=d_i, 1\leq i\leq n\}\,.
$$
Then it follows from (\ref{meas}) and (\ref{dist})
that
$$
\Pr\Big\{{\cal L}_n^{\hat{D}}=[k] \, \Big| \,\hat{D}_i=d_i, 1\leq i\leq n\Big\}=
\frac{\big|\{f\in{\cal M}_n(\vec{d}):{\cal L}_n(f)=[k]\}\big| }
{n!\big(\prod_{i=1}^n d_i!\big)^{-1}} .
$$
In order to count the set $\big\{f\in{\cal M}_n(\vec{d}) : {\cal L}_n(f)=[k]\big\}$,
we first construct a bijection between the set, ${\cal S}_{ \vec{d}}$,
of sequences ${\bf x}=(x_1, x_2,..., x_n)$ such that for any ${\bf x}\in{\cal S}_{\vec{d}}$
and $1\leq i\leq n$, we have $\big|\{m: x_m=i\}\big|=d_i$ and the
set  of mappings ${\cal M}_n(\vec{d})$. 

The bijection is defined in terms of an algorithm which,
for any sequence ${\bf x}=(x_1, x_2,..., x_n)\in {\cal S}_{\vec{d}}$,
constructs, after a series of `rounds', a corresponding mapping
 $f_{\bf x}\in{\cal M}(\vec{d})$. Informally,
the algorithm works as follows. In the first round the vertices
$i\in [n]$ with $d_i=0$ are mapped (in a way which is determined by the sequence ${\bf x}$) to  vertices in $[n]$ which have in-degree greater than 0.
After the first round the `availability' of some vertices to receive
directed edges in the corresponding directed graph
will have been reduced because the algorithm has mapped some vertices to them. A vertex $i$ becomes `unavailable' as soon as
it is `full', i.e. $d_i$ vertices have been mapped to it. So, at the beginning of the
second round some vertices which were `available' in the first round may be `unavailable' at the start of the second round. These new `unavailable'
vertices are mapped to `available' vertices in the second round and this
construction continues until the mapping $f_{\bf x}$ is
completely defined.

In order to describe rigorously the algorithm which constructs $f_{\bf x}$ , we introduce the following
notation. For any vertex $i\in [n]$, any $m\geq 1$,  and any ${\bf x}\in {\cal S}_{\vec{d}}$,
we let $a_i(m, {\bf x})$
denote the `availability' of vertex $i$ at the start of round $m$.
Initially we set $a_i(1, {\bf x})=d_i$ for all $i\in [n]$. The values of $a_i(m,{\bf x})$ when
$m\geq 2$ are determined recursively by the algorithm. For
${\bf x}\in {\cal S}_{\vec{d}}$ and $m\geq 1$ we define
the sets
$$
{\cal Y}_m({\bf x})=\{i\in [n]: a_i(m, {\bf x})=0\}, \quad {\cal Y}_0({\bf x})=\emptyset, \quad {\cal Z}_m({\bf x})= {\cal Y}_m({\bf x})\setminus {\cal Y}_{m-1}({\bf x}).$$
Finally, for ${\bf x}\in {\cal S}_{\vec{d}}$\,, let $y_0({\bf x})=0$ and  for $m\geq 1$, let $y_m({\bf x})=
|{\cal Y}_m({\bf x})|$ and $z_m({\bf x})=
|{\cal Z}_m({\bf x})|$.
\vskip .1in
\noindent {\bf The algorithm:} Given ${\bf x}\in{\cal S}_{\vec{d}}$\,, the corresponding mapping
$f_{\bf x}\in {\cal M}_n(\vec{d})$ is
constructed as follows.

\noindent {\bf Step 1.} Set $m=1$.
\vskip .05in
\noindent {\bf Step 2.} If ${\cal Z}_m({\bf x})\not=\emptyset$, then
\begin{itemize}
\item list the elements of ${\cal Z}_m({\bf x})$ in increasing order: $i_1<i_2<...<i_{z_m({\bf x})}$.
For $1\leq k\leq z_m({\bf x})$, define $f_{\bf x}(i_k)=x_{y_{m-1}({\bf x})+k}$,
where $x_t$ is the $t^{th}$ term in the sequence ${\bf x}=(x_1, x_2,...,x_n)$.
\item Next, for each $j\in [n]$, set
$$a_j(m+1, {\bf x})=a_j(m, {\bf x})-|\{i\in {\cal Z}_m({\bf x})
: f_{\bf x}(i)=j\}|.$$
\item Go to Step 3.
\end{itemize}

\noindent Otherwise, if ${\cal Z}_m({\bf x})=\emptyset$ (i.e.$ {\cal Y}_{m-1}({\bf x})=
{\cal Y}_{m}({\bf x})$), then
\begin{itemize}
\item list  the set
$[n]\setminus {\cal Y}_{m}({\bf x})$ in increasing order: $i_1<i_2<...<i_{n-y_{m}({\bf x})}$.
For $1\leq k\leq n-y_{m}({\bf x})$, let $f_{\bf x}(i_k)=x_{y_{m}({\bf x})+k}$. This completes
the construction of the mapping $f_{\bf x}$.
\end{itemize}

\vskip .05in
\noindent {\bf Step 3.} Set $m=m+1$ and go to Step 2.
\vskip .1in
\noindent We mention here that we have recently learned that Blitzstein and Diaconis \cite{BD06}
have used an algorithm  to construct labelled trees with a given degree sequence which 
is similar in spirit to the algorithm described above.
It is clear from the definition of the algorithm, that if ${\bf x},{\bf x}'\in {\cal S}_{\vec{d}}$
and ${\bf x}\not={\bf x}'$, then $f_{\bf x}\not=f_{{\bf x}'}$. In particular, if $k({\bf x}, {\bf x}')
\equiv \min\{k\in [n] : x_{k}\not= x_{k}'\}$ then there is some $i\in [n]$ such that
$f_{\bf x}(i)=x_{k({\bf x}, {\bf x}')}$ and $f_{\bf x'}(i)=x_{k({\bf x}, {\bf x}')}'$.
Since $\big|{\cal S}_{\vec{d}}\,\big|=\frac{n!}{d_1!\cdot\cdot\cdot d_n!}=\big|{\cal M}_n(\vec{d})\big|$,
it follows that
the algorithm constructs a bijection between ${\cal S}_{\vec{d}}$ and ${\cal M}_n(\vec{d})$.
Hence, to count $\{f\in{\cal M}_n(\vec{d}) : {\cal L}_n(f)=[k]\}$ it suffices to count
$\{{\bf x}\in {\cal S}_{\vec{d}}: {\cal L}_n(f_{\bf x})=[k]\}$. To this end, we prove the
following lemmas.
\begin{lem}
\label{Lemma1} For each ${\bf x}\in{\cal S}_{\vec{d}}$,
${\bf x}=(x_1, x_2, ..., x_n)$
we define $t({\bf x})$ as follows:
$$t({\bf x})=\min\big\{ t :
\big|\{x_{t}, x_{t+1},..., x_n\big\}\big|=n-t+1 \big\}\,.$$
Then
$$
{\cal L}_n\big(f_{\bf x}\big)=
\big\{x_{t({\bf x})}, x_{t({\bf x})+1}, ..., x_n\big\}.
$$
\end{lem}
\begin{proof}{\it (of Lemma 1)} Suppose that ${\bf x}\in{\cal S}_{\vec{d}}$ and
suppose that the algorithm terminates in round $\hat{m}$, then
the first step is to show that
\begin{equation}
\label{newlem1*}
 {\cal L}_n(f_{\bf x}) = [n]\setminus {\cal Y}_{\hat{m}}({\bf
x})\,.
\end{equation}
We begin with a few simple observations. First, it follows from the
description of the algorithm that at the beginning of the $m^{th}$
round, where $1\leq m\leq \hat{m}$, the vertices in $[n]\setminus {\cal
Y}_{m-1}({\bf x})$ have
not yet been assigned their image under the mapping $f_{\bf x}$.
Also, for every $i\in [n]$, we have
\begin{equation}
\label{lem2*} a_i(m, {\bf x})=\big|\{ k: x_k=i, y_{m-1}({\bf
x})+1\leq k\leq n\}\big|.
\end{equation}
In other words,  $a_i(m, {\bf x})=r$ if and only if $i$ appears $r$
times in the sequence $(x_{y_{m-1}({\bf x})+1}, x_{y_{m-1}({\bf
x})+2}, ..., x_n)$. In addition, it follows from (\ref{lem2*}) and
from the definition of ${\cal Y}_{m}({\bf x})$, that
$$
n-y_{m-1}({\bf x})=\sum_{i\in [n]} a_i(m, {\bf x})=\sum_{i\in [n]\setminus {\cal Y}_{m}({\bf x})} a_i(m, {\bf x})
$$
since $a_i(m, {\bf x})=0$ for every $i\in {\cal Y}_{m}({\bf x})$.
Now if $m=\hat{m}$, then  ${\cal
Y}_{\hat{m}}({\bf x})={\cal Y}_{\hat m-1}({\bf x})$ and so
$$
n-y_{\hat{m}}({\bf x})=n-y_{\hat{m}-1}({\bf x})=\sum_{i\in [n]\setminus {\cal Y}_{\hat{m}}({\bf x})} a_i(\hat{m},
{\bf x})\,.
$$
This implies immediately that
\begin{equation}
\label{newlem1+}
 a_i(\hat{m}, {\bf x})=1\,\mbox{ for all }
i\in [n]\setminus {\cal Y}_{\hat{m}}({\bf x}).
\end{equation}
 Hence the last step
of the algorithm generates a permutation on the set  $[n]\setminus
{\cal Y}_{\hat{m}}({\bf x})$  i.e.,
\begin{equation}
 \label{inclusion0}
\{x_{y_{\hat {m}({\bf x})+1}}, ..., x_{n}\}=[n]\setminus {\cal Y}_{\hat{m}}({\bf x}) \subseteq {\cal
L}_n(f_{\bf
x})\,.
\end{equation}
Suppose now that $i\in {\cal Y}_{\hat{m}}({\bf
x})$ (i.e., $i\in [n]$ and $a_i(m, {\bf x})=0$ for
some $1\leq m\leq \hat{m}-1$).  Then, since
$${\cal Y}_1({\bf x})\subseteq {\cal Y}_1({\bf x})
\subseteq ...\subseteq {\cal Y}_{\hat{m}-1}({\bf x})={\cal
Y}_{\hat{m}}({\bf x}),
$$
we can define
$$
m_i= \min\{ m: a_i(m, {\bf x})=0\}.
$$
 If $m_i=1$,
then $a_i(1,{\bf x})=d_i=0$, and we must have
$i\notin {\cal L}_n(f_{\bf x})$ since vertex $i$ has in-degree 0
in $G_n(f_{\bf x})$. So, suppose that $2\leq m_i\leq \hat{m}-1$ (and
$d_i=a_i(1, {\bf x})>0$).
It follows from the description of the
algorithm and the definition of $m_i$, that at the start of round $m_i$
vertex $i$ is
`unavailable' and $f_{\bf x }^{-1}(i)\equiv \{j: f_{\bf
x}(j)=i\}\subseteq {\cal Y}_{m_i-1}$.
Since the value of $f_{\bf x}(i)$ is determined during round $m_i$ and
since $a_i(m_i, {\bf x})=0$,
it follows that $f_{\bf x}(i)\not= i$.
 It is also clear from the
description of
the algorithm that for {\it any} $j\in [n]\setminus {\cal Y}_{m_i-1}$,
we  have $f_{\bf x}(j)\notin {\cal Y}_{m_i-1}$ since,
for every $k\in {\cal Y}_{m_i-1}$ and
for every $m\geq m_i$, ~$a_k(m, {\bf
x})=0$
(i.e., vertex $k\in {\cal Y}_{m_i-1}$ is `unavailable' in {\it every}
round $m\geq m_i$).
It follows by induction that
\begin{equation}
\label{lem1*}
 f_{\bf x}^{(t)}(i)\notin {\cal Y}_{m_i-1}
 \end{equation}
for every $t\geq 1$. Now suppose that $i\in {\cal L}_n(f_{\bf x})$,
and in particular, that $f_{\bf x}^{(t')}(i)=i$ for some $t'\geq 2$,
then $f^{(t'-1)}_{\bf x}(i)\in f^{-1}_{\bf x}(i)\subseteq {\cal
Y}_{m_i-1}$. But this contradicts (\ref{lem1*}), so we must have
$i\notin {\cal L}_n(f_{\bf x})$. Therefore
$$
{\cal Y}_{\hat{m}}({\bf x})\subseteq [n]\setminus  {\cal
L}_n(f_{\bf x})\,.
$$
This together with (\ref{inclusion0})  implies (\ref{newlem1*}).

\par
Next, we note that (\ref{newlem1+}), (\ref{inclusion0}), and the definition
of $t({\bf x})$, imply
\begin{equation}
\nonumber
 \{x_{y_{\hat {m}({\bf x})+1}}, ..., x_{n}\}=[n]\setminus {\cal Y}_{\hat{m}}({\bf x}) \subseteq \{x_{t({\bf x})}, ...., x_n\}.
\end{equation}
Thus to show that
\begin{equation}
\label{biglem*}
[n]\setminus {\cal Y}_{\hat m}({\bf x})=\{x_{t({\bf x})}, ...., x_n\},
\end{equation}
it suffices to show that
\begin{equation}
\nonumber x_{y_{\hat {m}({\bf x})}} \notin \{x_{t({\bf x})}, ...., x_n\}\,.
\end{equation}
Suppose that $x_{y_{\hat {m}({\bf x})}} \in \{x_{t({\bf x})}, ...., x_n\}$. Then in the ${\hat m}-1^{st}$
round there is some $j'\in {\cal Z}_{\hat{m}-1}$ such that $j'$ is assigned to
$x_{y_{\hat{m}({\bf x})}}$. So at the beginning of the ${\hat{m}-1}^{st}$
round of the algorithm  $x_{y_{\hat {m}({\bf x})}}$ was available, but it is unavailable at the start
of the ${\hat m}^{th}$ round (since $x_{y_{\hat {m}({\bf x})}}\notin [n]\setminus {\cal Y}_{\hat m}$).  
Since the
${\hat{m}-1}^{th}$ round of the algorithm cannot generate {\it new} unavailable vertices, we have a  contradiction. Hence (\ref{biglem*})
holds and  Lemma 1 follows immediately from (\ref{newlem1*}) .
\end{proof}

\noindent It follows from Lemma 1 that
$$\big\{{\bf x}\in{\cal S}_{\vec{d}}: {\cal L}_n(f_{\bf x})=[k]\big\}
=\big\{{\bf x}\in{\cal S}_{\vec{d}}: \{x_{n-k+1}, ..., x_n\}=[k]~\mbox{and}~x_{n-k}\in[k]\big\}.
$$
So
routine counting arguments yield
$$
\Pr\left\{{\cal L}_n^{\hat{D}}=[k] \bigg| \hat{D}_i=d_i, 1\leq i\leq n\right\}
$$
$$=
\frac{|\{{\bf x}\in{\cal S}_{\vec{d}}: \{x_{n-k+1}, ..., x_n\}=[k]~\mbox{and}~x_{n-k}\in[k]\}|\times\prod_{i=1}^n d_i!}
{n!}
$$
$$
={n\choose k}^{-1}\left({1\over n-k}\right)\sum_{i=1}^k d_1d_2\cdot\cdot\cdot d_{i-1}(d_i-1)d_i
d_{i+1}\cdot\cdot\cdot d_k
$$
for $1\leq k\leq n-1$.
Hence
$$
\Pr\big\{{\cal L}_n^{\hat{D}}=[k]\big\}=
{n\choose k}^{-1}\left({1\over n-k}\right)\sum_{i=1}^k E\big(\hat{D}_1\cdot\cdot\cdot
(\hat{D}_i-1)\hat{D}_i\cdot\cdot\cdot \hat{D}_k\big)
$$
$$
={n\choose k}^{-1}{k\over n-k}E\big((\hat{D}_1-1)\hat{D}_1\hat{D}_2\cdot\cdot\cdot \hat{D}_k\big)
$$
since $\hat{D}_1, \hat{D}_2, ..., \hat{D}_k$ are exchangeable. Repeating the
argument given above and using the exchangeability of the variables
$\hat{D}_1, ...., \hat{D}_n$, we also obtain
$$
\Pr\big\{{\cal L}_n^{\hat{D}}=L\big\}
={n\choose k}^{-1}{k\over n-k}E\big((\hat{D}_1-1)\hat{D}_1\hat{D}_2\cdot\cdot\cdot \hat{D}_k\big)
$$
for any $L\subseteq [n]$ such that $|L|=k$. Summing over all
such sets $L$,
we obtain
$$
\Pr\big\{X_n^{\hat{D}}=k\big\}=\frac{k}{n-k}E\big((\hat{D}_1-1)\hat{D}_1\hat{D}_2\cdot\cdot\cdot \hat{D}_k\big).
$$

Finally, we consider the case $k=n$. Clearly $X_n^{\hat{D}}=n$ if and only if
$f_{\bf x}$ is a permutation of $[n]$, and $f_{\bf x}$ is a permutation of $[n]$
if and only if $\hat{D}_i=1$ for all $1\leq i\leq n$. Hence
$$
\Pr\big\{X_n^{\hat{D}}=n\big\}= \Pr\big\{\hat{D}_i=1, 1\leq i\leq n\big\}
$$
as required, and the theorem is proved.
\end{proof}

Let $N_n^{\hat D}$ denote the number of components in $G_n^{\hat D}$.
Let $\sigma(m)$ be a uniform random permutation on an $m$-element set
and let $N_{\sigma(m)}$ denote the number of cycles in the
random permutation $\sigma(m)$. Then we have
\begin{cor}
\label{Corollary 1}
For $1\leq \ell\leq n$,
$$
\Pr\{N_n^{\hat D}=\ell\}
=\sum_{k=\ell}^n\frac{
|s(k,l)|}{(k-1)!(n-k)}E((\hat{D}_1-1)\hat{D}_1\hat{D}_2\cdot\cdot\cdot
\hat{D}_k)
$$
$$
\qquad \qquad+ \frac{|s(n,l)|}{n!}\Pr\{\hat{D}_i=1, 1\leq i\leq n\}\, ,
$$
where $s(\cdot \,,\cdot)$ are the Stirling numbers of the first kind.
\end{cor}

\begin{proof} Recall that the mapping $T_n^{\hat{D}}$
restricted to ${\cal L}_n^{\hat{D}}$, the cyclic vertices  of $T_n^{\hat{D}}$,
is a permutation of ${\cal L}_n^{\hat{D}}$. The corollary follows from the observation that
the number of components in $G_n^{\hat{D}}$ equals the number
of cycles in the permutation of ${\cal L}_n^{\hat{D}}$ by $T_n^{\hat{D}}$.
So $N_n^{\hat D}=\ell$
if and only if the mapping $T_n^{\hat D}$ restricted to ${\cal L}_n^{\hat{D}}$
is a permutation with $\ell$ cycles.
\par
First suppose that $1\leq \ell \leq k\leq n-1$, and $L\subseteq [n]$ such that $|L|=k$. Also fix the
in-degree
sequence $\vec{d}=(d_1, d_2,..., d_n)$ such that $\sum_{i=1}^n d_i=n$,  and such that
$d_i\geq 1$ for every $i\in L$ and $d_{i'}\geq 2$ for {\it some} $i'\in L$.
Recall that
for any $f,g\in{\cal M}_n(\vec{d})$ we have ${\bf P}_n^{\hat{D}}(f)={\bf P}_n^{\hat{D}}(g)$,
so it follows from the bijection described in the proof of Theorem 1 that
$$\Pr\Big\{N^{\hat D}_n=\ell \, \Big|\, {\cal L}_n^{\hat D}=L, \hat{D}_i=d_i, 1\leq i\leq n\Big\}=$$
$$={\Big|\big\{{\bf x}\in {\cal S}_{\vec d} : \{ x_{n-k+1}, ..., x_n\}=L, x_{n-k}\in L,
\mbox{permutation }f_{\bf x}\big|_L
\mbox{ has }\ell\mbox{ cycles}\big\}\Big|\over \Big|\big\{{\bf x}\in {\cal S}_{\vec
d} : \{ x_{n-k+1}, ..., x_n\}=L, x_{n-k}\in L\big\}\Big|}$$
$$=\Pr\{N_{\sigma(k)}=\ell\}.$$
Since the set $L\subseteq [n]$ and the degree sequence $\vec{d}$ were
arbitrary, it follows that
\begin{equation}
\label{cor1_eq1}
\Pr\left\{N_n^{\hat D}=\ell \big| X_n^{\hat D}=k\right\} =\Pr\{N_{\sigma(k)}=\ell\}.
\end{equation}
In the case  $X_n^{\hat D}=n$, the mapping $T_n^{\hat D}$ is a
random permutation of $[n]$ (since ${\cal L}_n^{\hat D}=[n]$), so for
$1\leq \ell\leq n$ we also have
\begin{equation}
\label{cor1_eq2}
\Pr\left\{N_n^{\hat D}=\ell \big| X_n^{\hat D}=n\right\} =\Pr\{N_{\sigma(n)}=\ell\}.
\end{equation}
It follows that
$$
\Pr\{N_n^{\hat D}=\ell\}=\sum_{k=\ell}^n
\Pr\left\{N_n^{\hat D}=\ell \big| X_n^{\hat D}=k\right\}\Pr\{X_n^{\hat D}=k\}
$$
$$
=\sum_{k=1}^n \Pr\{N_{\sigma(k)}=\ell\}\Pr\{X_n^{\hat{D}}=k\}\,.$$
Let $s(\cdot \,,\cdot)$ denote the Stirling numbers of the first kind, then it is well known that there are
$|s(k,l)|$ permutations of $k$-element set with exactly $l$ cycles, i.e.,
$$
\Pr\{N_{\sigma(k)}=\ell\} = \frac{
|s(k,l)|}{k!}\,,
$$
which implies the assertion of the corollary.
\end{proof}
\par
Let ${\cal B}_n^{\hat{D}}$ denote the event that the random graph
$G_n^{\hat{D}}$ is connected. Then since ${\cal B}_n^{\hat D}=\{N_n^{\hat D}=1\}$, we obtain 
the following result immediately
 from the proof of Corollary 1:

\begin{cor}
\label{Corollary2}
$$
\Pr\big\{{\cal B}_n^{\hat{D}}\big\}=\sum_{k=1}^n \frac{1}{k}\Pr\big\{X_n^{\hat{D}}=k\big\}
$$
$$
=\sum_{k=1}^{n-1}{E\big((\hat{D}_1-1)\hat{D}_1\hat{D}_2\cdot\cdot\cdot \hat{D}_k\big)\over n-k}
+{\Pr\big\{\hat{D}_i=1,1\leq i\leq n\big\}\over n}.
$$
\end{cor}
Finally, suppose that $\xi_1, \xi_2,...$ is a sequence of independent indicator
variables such that, for $k\geq 1$, $\Pr\{\xi_k=1\}={1\over k}$ and such that
$\xi_1, \xi_2, ....$ and $X_n^{\hat D}$ are independent. It is well known
(see \cite{Feller70}) that for $m\geq 1$
\begin{equation}
\label{sigma}
N_{\sigma(m)}\stackrel{d}{ \sim} \sum_{k=1}^m \xi_k,
\end{equation}
and  that $(N_{\sigma(m)}-\log m)/\sqrt{\log m}$ converges in distribution to the
standard $N(0,1)$ distribution.
It is an easy consequence of (\ref{cor1_eq1})-(\ref{sigma}), that
\begin{cor}
\label{Corollary 3}
For $n\geq 1$,
\begin{equation}
\label{numbercomp}
N_n^{\hat D}\stackrel{d}{ \sim}\sum_{k=1}^{X_n^{\hat D}}\xi_k.
\end{equation}
\end{cor}
\noindent We note that (\ref{numbercomp}) is useful for investigating
the asymptotic distribution of $N_n^{\hat D}$.
\par Next, we consider the distribution of the size of a `typical'
component of $G_n^{\hat D}$. Let $C_1^{\hat D}(n)$ denote the size of the component in
$G_n^{\hat D}$ which contains the vertex 1, then the distribution of
$C_1^{\hat D}(n)$ is given by the following theorem.

\begin{thm}
\label{Theorem2}
For $1\leq \ell\leq n$, let $D_1', D_2', ..., D_{\ell}'$ be a sequence of variables
with joint distribution given by
$$\Pr\Big\{ D_i'=d_i, ~1\leq i\leq \ell\Big\}=\Pr\Big\{ {\hat D}_i=d_i,   ~ 1\leq i\leq \ell \, \Big| \,
\sum_{i=1}^{\ell} {\hat D}_i=\ell\Big\},
$$
then
$$
\Pr\Big\{C_1^{\hat D}(n)=\ell\Big\}=\frac{\ell}{n}\,\Pr\Big\{{\cal B}^{D'}_{\ell}\Big\}
\,\Pr\Big\{\sum_{i=1}^{\ell} {\hat D}_i=\ell\Big\} .
$$
\end{thm}
\begin{proof} Fix $1\leq \ell\leq n$ and
let ${\cal C}_1^{\hat D}(n)$  denote the vertex set of the component of $G_n^{\hat D}$
which contains the vertex $1$.
Then we have
\begin{equation*}
\Pr\{C_1^{\hat D}(n)=\ell\}=\sum_{C\subseteq [n] ~s.t.~ 1\in C\atop and ~ |C|=\ell}
\Pr\{{\cal C}_1^{\hat D}(n)=C\}.
\end{equation*}
Now fix $C=[\ell]$ and observe that if ${\cal C}_1^{\hat D}(n)=[\ell]$
then $T_n^{\hat D}$ must map $[\ell]$ into $[\ell]$,  and,
in particular, we must have
$\sum_{i=1}^{\ell} {\hat D}_i=\ell$. So
\begin{equation}
\label{thm2+}
\Pr\{{\cal C}_1^{\hat D}(n)=[\ell]\}=\Pr\Big\{{\cal C}_1^{\hat D}(n)= [\ell] \, \Big|\,
\sum_{i=1}^{\ell} {\hat D}_i=\ell\Big\}\Pr\Big\{\sum_{i=1}^{\ell}{\hat D}_i=\ell\Big\}.
\end{equation}
Next, by summing over all degree sequences $\vec{d}=(d_1, d_2, ..., d_n)$ such that $\sum_{i=1}^{\ell} d_i=\ell$
and $\sum_{i=1}^nd_i=n$, we obtain
$$
\Pr\Big\{{\cal C}_1^{\hat D}(n)= [\ell] \, \Big|\,
\sum_{i=1}^{\ell} {\hat D}_i=\ell\Big\}
$$
$$
=\sum_{\vec{d}~ s.t. \sum_{i=1}^{\ell} d_i=\ell,\atop
\sum_{i=1}^nd_i=n}
\Pr\Big\{{\cal C}_1^{\hat D}(n)=[\ell] \, \Big| \,  {\hat D}_i=d_i, 1\leq i\leq n\Big\}\quad\quad
$$
\begin{equation}
\label{thm2*}
\quad\quad\quad\quad\quad\times
\Pr\Big\{{\hat D}_i=d_i, 1\leq i\leq n \, \Big| \, \sum_{i=i}^{\ell}{\hat D}_i=\ell
\Big\}.
\end{equation}
Now fix $\vec{d}=(d_1, ..., d_n)$ such that $\sum_{i=1}^{\ell}d_i=\ell$ and
$\sum_{i=1}^n d_i=n$, and define $M(d_1, d_2, ..., d_{\ell})$ to be the number of mappings
$f:[\ell]\to [\ell]$ such that  $d_i(f)=d_i$ for $1\leq i\leq \ell$ and $G_{\ell}(f)$ is connected. Then
we have
$$
\Pr\Big\{{\cal C}_1^{\hat D}(n)=[\ell]\, \Big| \,  {\hat D}_i=d_i, 1\leq i\leq n\Big\}
={M(d_1, ..., d_{\ell})(n-\ell)!\over d_{\ell+1}!
\cdot\cdot\cdot d_n!}\times {d_1!\cdot\cdot\cdot d_n!\over n!}
$$
$$
={n\choose \ell}^{-1}M(d_1, ..., d_{\ell})\times{d_1!\cdot\cdot\cdot d_{\ell}!\over
\ell!}\quad\quad
$$
\begin{equation}
\label{thm2**}
={n\choose \ell}^{-1}
\Pr\Big\{{\cal B}^{{\cal D}'}_{\ell} \, \Big|\,  D'_i=d_i, 1\leq i\leq \ell\Big\}.
\end{equation}
Now for any degree sequence $\vec{d}=(d_1, ..., d_n)$, we define
$\vec{d}_{\ell}=(d_1, ..., d_{\ell})$ and ${\vec{d}_{\ell}}^{~\prime}=(d_{\ell+1}, ..., d_n)$
(so $\vec{d}=(\vec{d}_{\ell}, \vec{d}_{\ell}^{~\prime})$). Then, by substituting the
RHS of (\ref{thm2**}) into (\ref{thm2*}),
we obtain
$$
\Pr\Big\{{\cal C}_1^{\hat D}(n)= [\ell] \, \Big| \,
\sum_{i=1}^{\ell} {\hat D}_i=\ell\Big\}
$$
$$
={n\choose \ell}^{-1}\sum_{\vec{d}_{\ell}~ s.t. \atop\sum_{i=1}^{\ell} d_i=\ell}
\sum_{\vec{d}_{\ell}^{~\prime}~s.t.\atop \sum_{i=\ell+1}^n d_i=n-\ell}
\Pr\Big\{{\cal B}^{{\cal D}'}_{\ell} \, \Big| \,  D'_i=d_i, 1\leq i\leq \ell\Big\}
$$
$$\quad\quad\quad\quad\quad\times
\Pr\Big\{{\hat D}_i=d_i, 1\leq i\leq n \, \Big|\, \sum_{i=i}^{\ell}{\hat D}_i=\ell
\Big\}
$$
$$
={n\choose \ell}^{-1}\sum_{\vec{d}_{\ell}~ s.t. \atop\sum_{i=1}^{\ell} d_i=\ell}
\Pr\Big\{{\cal B}^{{\cal D}'}_{\ell} \, \Big|\,  D'_i=d_i, 1\leq i\leq \ell\Big\}
\Pr\Big\{{ D}_i'=d_i, 1\leq i\leq \ell
\Big\}
$$
\begin{equation}
\label{thm2***}
={n\choose \ell}^{-1}
\Pr\{{\cal B}^{{\cal D}'}_{\ell} \}.
\end{equation}
Substituting  (\ref{thm2***}) into (\ref{thm2+}), we obtain
$$
\Pr\Big\{{\cal C}_1^{\hat D}=[\ell]\Big\}={n\choose \ell}^{-1}
\Pr\Big\{{\cal B}^{{\cal D}'}_{\ell} \Big\}\Pr\Big\{\sum_{i=1}^{\ell}{\hat D}_i=\ell\Big\}.
$$
It follows by the same argument and from the exchangeability of the variables
${\hat D}_1, {\hat D}_2, ..., {\hat D}_n$ that, for any $C\subseteq [n]$ such that
$1\in C$ and $|C|=\ell$, we have
$$
\Pr\Big\{{\cal C}_1^{\hat D}=C\Big\}={n\choose \ell}^{-1}
\Pr\Big\{{\cal B}^{{\cal D}'}_{\ell} \Big\}\Pr\Big\{\sum_{i=1}^{\ell}{\hat D}_i=\ell\Big\}.
$$
So
$$\Pr\Big\{C_1^{\hat D}(n)=\ell\Big\}={n-1\choose \ell-1}{n\choose \ell}^{-1}
\Pr\Big\{{\cal B}^{{\cal D}'}_{\ell} \Big\}\Pr\Big\{\sum_{i=1}^{\ell}{\hat D}_i=\ell\Big\}
$$
$$
={\ell\over n}\Pr\Big\{{\cal B}^{{\cal D}'}_{\ell}\Big \}\Pr\Big\{\sum_{i=1}^{\ell}{\hat D}_i=\ell\Big\}
$$
as desired.
\end{proof}
\noindent We prove an extension of Theorem 2 in the following special case:  Suppose that $D_1, D_2, ...$
is a sequence of i.i.d. non-negative integer valued random variables. For each $n\geq 1$,
let $\hat D(n)=(\hat D_{1,n}, ...., \hat D_{n,n})$ be a sequence of variables with joint distribution
given by
$$
\Pr\Big\{\hat D_{i,n}=d_i\mbox{ for } i=1,2,.., n\Big\}=
\Pr\Big\{ D_i=d_i,~~i=1,2,...,n\, \Big| \,  \sum_{i=1}^n D_i =n\Big\}.
$$
Let ${\cal C}_1^{\hat D(n)}$ denote the vertex set of the connected component
in $G_n^{\hat D(n)}\equiv G(T_n^{\hat D(n)})$ which contains the vertex labelled
$1$. For $k>1$, we define ${\cal C}_k^{\hat D(n)}$ recursively as follows:
If $[n]\setminus ({\cal C}_1^{\hat D(n)}\cup\cdot\cdot\cdot\cup {\cal C}_{k-1}^{\hat D(n)})
\not=\emptyset$, let
${\cal C}_k^{\hat D(n)}$ denote the vertex set of the connected component in $G_n^{\hat D(n)}$
which contains the smallest element of $[n]\setminus ({\cal C}_1^{\hat D(n)}\cup\cdot\cdot\cdot\cup {\cal C}_{k-1}^{\hat D(n)})$;  otherwise,
set ${\cal C}_k^{\hat D(n)}=\emptyset$. For all $k\geq 1$, let $C_k^{\hat D(n)}=|{\cal C}_k^{\hat D(n)}|$,
then we have
\begin{thm}
\label{Theorem3}
Suppose $ 1\leq k< n$ and  
$\ell_1, \ell_2, \ldots , \ell_k$ are
such that
$1\leq \ell_1<n$\,,
$1\leq \ell_i<n-\ell_0-\ell_1-\cdot\cdot\cdot -\ell_{i-1}\leq n$,
for $i=2,..., k$, and $\sum_{i=1}^k\ell_i\leq n$. Then we have
$$
\Pr\big\{ C_1^{\hat D(n)}=\ell_1, ..., C_k^{\hat D(n)}=\ell_k \big\}=
$$
$$
=\left[\prod_{i=1}^k\left({\ell_i\over n-t_{i-1}}\right)
\Pr\left\{{\cal B}_{\ell_i}^{\hat D(\ell_i)}\right\}\right]\times
\Pr\left\{\sum_{j=t_{i-1}+1}^{t_i}\hat D_{j,n}=\ell_i,~1\leq i\leq k\right\},
$$
where $t_0=0$ and $t_i\equiv \ell_1+...+\ell_i$ , \,  $i=1,2,..., k$.
\end{thm}
\begin{proof} For any integers $0<m_1<m_2$, let  $[m_1, m_2]$ denote the
set of integers $\{m_1, m_1+1, ..., m_2\}$. Then we have
$$
\Pr\{ C_1^{\hat D(n)}=\ell_1, ..., C_k^{\hat D(n)}=\ell_k\}=
$$
\begin{equation}
\label{thm3_eq1}
=\prod_{i=1}^k{n-t_{i-1}-1\choose \ell_i-1}\times
\Pr\left\{{\cal C}_i^{\hat D(n)}=[ t_{i-1}+1, t_i]\mbox{ for }i=1,2,...,k\right\}.
\end{equation}
Since
$$\left\{{\cal C}_i^{\hat D(n)}=[ t_{i-1}+1, t_i],~1\leq i\leq k\right\}\subseteq
\left\{\sum_{j=t_{i-1}+1}^{t_i}\hat D_{j,n}=\ell_i,~1\leq i\leq k\right\}
$$
we have
$$
\Pr\left\{{\cal C}_i^{\hat D(n)}=[ t_{i-1}+1, t_i]\mbox{ for }i=1,2,...,k\right\}
$$
$$
=
\Pr\left\{{\cal C}_i^{\hat D(n)}=[ t_{i-1}+1, t_i]\mbox{ for }i=1,2,...,k
\bigg|\sum_{j=t_{i-1}+1}^{t_i}\hat D_{j,n}=\ell_i,~1\leq i\leq k \right\}
$$
\begin{equation}
\label{thm3_eq2}
\times
\Pr\left\{\sum_{j=t_{i-1}+1}^{t_i}\hat D_{j,n}=\ell_i,~1\leq i\leq k\right\}.
\end{equation}
Let ${\cal A}=\{\vec{d}=(d_1, .., d_n) : \sum_{j=t_{i-1}+1}^{t_i}d_j=\ell_i, 1\leq i\leq k
\mbox{ and }\sum_{j=1}^nd_j=n\}$, then
$$
\Pr\left\{{\cal C}_i^{\hat D(n)}=[ t_{i-1}+1, t_i]\mbox{ for }i=1,2,...,k
\bigg|\sum_{j=t_{i-1}+1}^{t_i}\hat D_{j,n}=\ell_i,~1\leq i\leq k \right\}
$$
$$
=\sum_{\vec{d}\in{\cal A}}\Pr\left\{{\cal C}_i^{\hat D(n)}=[ t_{i-1}+1, t_i]\mbox{ for }i=1,2,...,k\bigg|
\hat D_{j,n}=d_j, 1\leq j\leq n\right\}
$$
\begin{equation}
\label{thm3_eq3}
\quad\quad\quad\times
\Pr\left\{\hat D_{j,n}=d_j, 1\leq j\leq n\bigg |\sum_{j=t_{i-1}+1}^{t_i}\hat D_{j,n}=\ell_i,~1\leq i\leq k
\right \}\,.
\end{equation}
Now for any $1\leq i\leq k$, let $M(d_{t_{i-1}+1}, d_{t_{i-1}+2},...., d_{t_i})$ denote the number of
mappings of $[t_{i-1}+1, t_i]$ into $[t_{i-1}+1, t_i]$ such that the corresponding directed graph
on the vertex set $\{t_{i-1}+1,..., t_i\}$ is connected and the in-degree of each vertex
$t_{i-1}+1\leq m\leq t_i$ is $d_m$. Then for any $\vec{d}\in{\cal A}$, we have
(as in the proof of Theorem 2)
$$
\Pr\left\{{\cal C}_i^{\hat D(n)}=[ t_{i-1}+1, t_i]\mbox{ for }i=1,2,...,k\bigg|
\hat D_{j,n}=d_j, 1\leq j\leq n\right\}
$$
$$=\prod_{i=1}^k{n-t_{i-1}\choose \ell_i}^{-1}\times
\prod_{i=1}^k{M(d_{t_{i-1}+1},...,d_{t_i})(d_{t_{i-1}+1})!\cdot\cdot\cdot d_{t_i}!\over
\ell_i!}
$$
\begin{equation}
\label{thm3_eq4}
=\prod_{i=1}^k{n-t_{i-1}\choose \ell_i}^{-1}\times
\prod_{i=1}^k\Pr\left\{{\cal B}_{\ell_i}^{\hat D(\ell_i)} \bigg|
\hat D_{j,\ell_i}=d_{t_{i-1}+j}\mbox{ for }1\leq j\leq \ell_i\right\}.
\end{equation}
The last equality follows by `re-labelling' , for each $1\leq i\leq k$,
the vertices $t_{i-1}+1,..., t_{i}$ by $1,2,..., \ell_i$.

Next  we note that since the variables $D_1, D_2, ..., D_n$ are independent and identically
distributed, we have for every $\vec{d}\in{\cal A}$
$$
\Pr\left\{\hat D_{j,n}=d_j, 1\leq j\leq n\bigg |\sum_{j=t_{i-1}+1}^{t_i}\hat D_{j,n}=\ell_i,~1\leq i\leq k
\right \}
$$
$$
=\prod_{i=1}^{k}{\Pr\{D_j=d_j, t_{i-1}+1\leq j\leq t_i\}\over \Pr\{\sum_{j=t_{i-1}+1}^{t_i}D_j=\ell_i\}}
\times{\Pr\{D_j=d_j,~ t_k+1\leq j\leq n\}\over \Pr\{\sum_{j=t_{k}+1}^nD_j= n-t_k\}}
$$
\begin{equation}
\label{thm3_eq5}
=\prod_{i=1}^{k}\Pr\left\{\hat D_{j,\ell_i}=d_{t_{i-1}+j}, ~1\leq j\leq \ell_i\right\}
\times{\Pr\{D_j=d_j,~ t_k+1\leq j\leq n\}\over \Pr\{\sum_{j=t_{k}+1}^nD_j= n-t_k\}}
\end{equation}
where we replace the the last factor in the product by $1$ if $t_k=n$. Substituting (\ref{thm3_eq4})
and (\ref{thm3_eq5}) into (\ref{thm3_eq3}) and summing over all $\vec{d}\in {\cal A}$,
we obtain
$$
\Pr\left\{{\cal C}_i^{\hat D(n)}=[ t_{i-1}+1, t_i]\mbox{ for }i=1,2,...,k
\bigg|\sum_{j=t_{i-1}+1}^{t_i}\hat D_{j,n}=\ell_i,~1\leq i\leq k \right\}
$$
\begin{equation}
\label{thm3_eq6}
=\prod_{i=1}^k{n-t_{i-1}\choose \ell_i}^{-1}\times
\prod_{i=1}^k\Pr\left\{{\cal B}_{\ell_i}^{\hat D(\ell_i)}\right\}.
\end{equation}
Finally, substituting (\ref{thm3_eq6}) into (\ref{thm3_eq2}) and then substituting
(\ref{thm3_eq2}) into (\ref{thm3_eq1}), we obtain the result.
\end{proof}

\par
In the general case where the
variables $\hat D_1, \hat D_2, ..., \hat D_n$ are exchangeable and $\sum_{i=1}^n
\hat D_i=n$, an analogue of Theorem 3 can be proved. However, the result for the special
case of Theorem 3 is easier because we are able to exploit the fact
that the variables $D_1, D_2, ..$ are i.i.d. to obtain the product in (\ref{thm3_eq5})
and hence the product in (\ref{thm3_eq6}).
Theorems 1,2, and 3 and their corollaries illustrate how the distributions of
random mapping statistics for $T_n^{\hat D}$ can be computed in terms
of the variables $\hat{D}_1, \hat{D}_2, ..., \hat{D}_n$. Similar results for
local properties of $T_n^{\hat D}$ such as the number of predecessors
and the number of successors of a given vertex (or vertices) are obtained
in a companion paper \cite{HJ06a}.


\section{Examples}
\label {s:examples}
We consider two special examples which correspond, respectively, to a random mapping with
`preferential attachment' and a random mapping with `anti-preferential
attachment'.
\vskip .1in
\subsection{\bf  A Preferential Attachment Model}
\vskip .1in
\noindent In this section we investigate $T_n^{\rho}:[n]\to[n]$, a random mapping with `preferential
attachment', where $\rho>0$ is a fixed parameter.
 For $1\leq k\leq n$, we define $T^{\rho}_n(k)=X_k^{(\rho,n)}$ where
 $X_1^{(\rho, n)}, X_2^{(\rho, n)},..., X_n^{(\rho,n)}$ is a sequence of random variables whose
 distributions depend on the evolution of an urn scheme.
The distribution of each $X_k^{({\rho},n)}$ is determined by a
(random) $n$-tuple of non-negative weights $\vec{a}(k)=(a_1(k), a_2(k), ..., a_n(k))$
where, for $1\leq j\leq n$, $a_j(k)$ is the `weight' of the $j^{th}$ urn
at the {\it start} of the $k^{th}$ round of the urn scheme.
Specifically, given $\vec{a}(k)=\vec{a}=(a_1, ..., a_n)$, we define
$$
\Pr\left\{X_k^{({\rho},n)}=j \big| \vec{a}(k)=\vec{a}\right\}={a_j\over \sum_{i=1}^{n}a_i}.
$$
The random weight vectors $\vec{a}(1), \vec{a}(2), ..., \vec{a}(n)$ associated with the
urn scheme are determined
recursively. For $k=1$, we set $a_1(1)=a_2(1)=\cdot\cdot\cdot=a_n(1)={\rho}>0$.
For $k>1$, $\vec{a}(k)$ depends on both $\vec{a}(k-1)$ and the value
of $X_{k-1}^{({\rho},n)}$ as follows: Given that $X_{k-1}^{({\rho},n)}=j$, we set $a_j(k)=a_j(k-1)+1$
and for all other $i\not=j$, we set $a_i(k)=a_i(k-1)$ (i.e. if $X_{k-1}^{(\rho,n)}=j$
then a `ball, with weight 1 is added to the $j^{th}$ urn).
\par
The random mapping $T_n^{\rho}$ as defined above is a preferential attachment model
in the following sense. Since, for $1\leq k\leq n$, $T_n^{\rho}(k)=X_k^{({\rho},n)}$,
and
since the (conditional) distribution of $X_k^{({\rho},n)}$ depends
on the state of the urn scheme at the start of round $k$,
it is clear that vertex $k$ is more likely to be mapped to vertex $j$ if the weight $a_j(k)$
is (relatively) large, i.e. if several of the vertices $1,2, ..., k-1$
have already been mapped to vertex $j$ . In the following proposition we establish the distribution
of $T_n^{\rho}$.

\begin{thm}
\label{Theorem 4}
Suppose that $D_1^{\rho}, D_2^{\rho}, ...$ are i.i.d.
random
variables with
a generalized negative binomial
distribution given by
$$
\Pr\{D_1^{\rho}=k\}={\Gamma(k+\rho)\over k!\Gamma(\rho)}
\left({\rho\over 1+\rho}\right)^{\rho}\left({1\over 1+\rho}\right)^k\quad\mbox{for}\quad k=0, 1, ...,
$$
where $\rho>0$ is a fixed parameter.

\noindent
For $n\geq 1$, let ${{\hat D}}({\rho},n)=({\hat D}^{\rho}_{1,n}, {\hat D}^{\rho}_{2,n}, ...,
{\hat D^{\rho}}_{n,n})$ be a sequence of variables  with joint distribution given by
$$
\Pr\{{\hat D}^{\rho}_{i,n}=d_i, 1\leq i\leq n\}=\Pr\left\{ D^{\rho}_i=d_i, 1\leq i\leq n\bigg|
\sum_{i=1}^n D^{\rho}_i= n\right\}.
$$
 Then for every $n\geq 1$, the random mappings $T_n^{\rho}:[n]\to [n]$ and
 $T_n^{\hat D(\rho,n)}:[n]\to [n]$ have the same distribution.

\end{thm}
\begin{proof}
To prove the result it is enough to show that for any $n\geq 1$ and any  $f\in{\cal M}_n$
$$
\Pr\{T_n^{\rho}=f\}=\Pr\{T_n^{\hat D(\rho, n)} =f\}.
$$
Suppose that $f\in{\cal M}_n$ and that $\vec{d}(f)=(d_1, d_2, ..., d_n)$.
It is straightforward to check that
$$
\Pr\{{\hat D}^{\rho}_{i,n}=d_i, 1\leq i\leq n\}=\Pr\left\{ D^{\rho}_i=d_i, 1\leq i\leq n\bigg|
\sum_{i=1}^n D^{\rho}_i= n\right\}
$$
$$
=\prod_{i=1}^n{\Gamma(d_i+\rho)\over (d_i)!\Gamma(\rho)}
\times {n!\Gamma(n\rho)\over \Gamma(n+n\rho)}
$$
since the distribution of
$\sum_{i=1}^n D^{\rho}$
is given by
\begin{equation}
\label{prop1}
\Pr\left\{\sum_{i=1}^nD^{\rho}_i=k\right\}={\Gamma(k+n\rho)\over k!\Gamma(n\rho)}
\left({\rho\over 1+\rho}\right)^{n\rho}\left({1\over 1+\rho}\right)^k\quad\mbox{for}\quad k=0, 1, ...\quad .
\end{equation}
It follows from the definition of $T_n^{\hat D(\rho, n)}$ that
$$
\Pr\{T_n^{\hat D(\rho, n)}=f\} = {d_1!\cdot\cdot\cdot d_n!\over n!}\times
\prod_{i=1}^n{\Gamma(d_i+\rho)\over (d_i)!\Gamma(\rho)}
\times {n!\Gamma(n\rho)\over \Gamma(n+n\rho)}$$
$$
={\prod_{i=1}^n(d_i+\rho-1)_{d_i}\over (n\rho+n-1)_n}
$$
where
$(m)_k\equiv m(m-1)\cdot\cdot\cdot(m-k+1)$ (and $(m)_0=1$).
\par
Next, observe that  if $T_n^{\rho}=f$  and $\vec{d}(f)=(d_1, d_2, ..., d_n)$, then
for each $1\leq i\leq n$,  $d_i$ `balls' of weight 1
are added to  $i^{th}$ urn during the evolution of the urn scheme described above.
So it follows from the definition of $T_n^{\rho}$ in terms of the
urn scheme,
that
$$
\Pr\{T_n^{\rho}=f\}={\prod_{i=1}^n(d_i+\rho-1)_{d_i}\over (n\rho+n-1)_n}.
$$
The result follows since $n\geq 1$ and $f\in{\cal M}_n$ were arbitrary.
\end{proof}
It is clear from Theorem 4 that the order in which a realisation of
$T_n^{\rho}$ is sequentially constructed does not matter.
In particular, suppose that $i_1, i_2, ..., i_n$ is a permutation of $[n]$
and for $1\leq k\leq n$, let $\tilde{T}_n^{\rho}(i_k)\equiv X_n^{(\rho, n)}$
where the variables $X_1^{(\rho, n)}, X_2^{(\rho, n)}, ..., X_n^{(\rho, n)}$
are as defined above. Then it follows from the proof of Theorem 4
that $\tilde{T}_n^{\rho}\stackrel{d}{ \sim} T_n^{\hat{D}(\rho,n)}\stackrel{d}{ \sim} T_n^{\rho}$.

\par Since
 $T_n^{\rho}\stackrel{d}{ \sim} T_n^{\hat D({\rho},n)}$, we can investigate the structure
 of $G_n^{\rho}\equiv G(T_n^{\rho})$ by considering the structure
 of $G_n^{\hat{D}(\rho, n)}$. In this paper, in order to illustrate the general method,
 we derive both exact and asymptotic results for
 $X_n^{\rho}$, the number of cyclic vertices in $G_n^{\rho}$,
 and $C_1^{\rho}(n)$, the size of the component in
 $G_n^{\rho}$ which contains vertex 1. Further results concerning the
 structure of $G_n^{\rho}$ can be found in \cite{HJ06b} where we also
 consider the structure of $G_n^{\rho}$ when $\rho=\rho(n)$ is a function
 of $n$.

 By Theorem 1 and Theorem 4 we have, for $1\leq k<n$,
\begin{equation}
\label{cyclic_1}
\Pr\Big\{X_n^{\rho}=k\}=\Pr\{X_n^{\hat D({\rho},n)}=k\Big\}=
{k\over n-k}E\Big((\hat D^{\rho}_{1,n}-1)\hat D^{\rho}_{1,n}\hat D^{\rho}_{2,n}\cdot
\cdot\cdot \hat D^{\rho}_{k,n}\Big).
\end{equation}
Since
$$
E\Big((\hat D^{\rho}_{1,n}-1)\hat D^{\rho}_{1,n}\hat D^{\rho}_{2,n}\cdot
\cdot\cdot \hat D^{\rho}_{k,n}\Big)=E\Big(( D^{\rho}_{1}-1) D^{\rho}_{1} D^{\rho}_{2}\cdot
\cdot\cdot  D^{\rho}_{k}\, \Big|\,  \sum_{i=1}^nD_i^{\rho}=n\Big)
$$
we have
$$
E\Big((\hat D^{\rho}_{1,n}-1)\hat D^{\rho}_{1,n}\hat D^{\rho}_{2,n}\cdot
\cdot\cdot \hat D^{\rho}_{k,n}\Big)=
$$
$$={\mbox{coefficient of }s^n\mbox{ in }
E\Big(( D^{\rho}_{1}-1) D^{\rho}_{1}s^{D_1^{\rho}} D^{\rho}_{2}s^{D_2^{\rho}}\cdot
\cdot\cdot  D^{\rho}_{k}s^{D_k^{\rho}}s^{D_{k+1}^{\rho}}\cdot\cdot\cdot
s^{D_n^{\rho}}\Big)\over
\mbox{coefficient of }s^n\mbox{ in }
E\Big(s^{D_{1}^{\rho}}\cdot\cdot\cdot
s^{D_n^{\rho}}\Big)}
$$
\begin{equation}
\label{cyclic_2}
={\mbox{coefficient of }s^n\mbox{ in }
E\Big(( D^{\rho}_{1}-1) D^{\rho}_{1}s^{D_1^{\rho}}\Big)\Big(E\big( D^{\rho}_{1}s^{D_1^{\rho}}\big)\Big)^{k-1}
\Big(E\big(s^{D_{1}^{\rho}}\big)\Big)^{n-k}\over
\mbox{coefficient of }s^n\mbox{ in }
\Big(E\big(s^{D_{1}^{\rho}}\big)\Big)^{n}}.
\end{equation}
The last equality holds since the variables $D_1^{\rho}, D_2^{\rho}, ..., D_n^{\rho}$
are independent and identically distributed. Since the identity
$$
{1\over (1-u)^{\alpha}}=\sum_{k=0}^{\infty}{\Gamma(k+\alpha)\over k!\Gamma(\alpha)}u^k
$$
holds for all $\alpha>0$ and $|u|<1$, we have
$$
E\Big(s^{D_1^{\rho}}\Big)=\left({\rho\over 1+\rho -s}\right)^{\rho},\quad
E\Big(D_1^{\rho}s^{D_1^{\rho}}\Big)=s\left({\rho\over 1+\rho -s}\right)^{\rho+1},\quad \mbox{ and}
$$
\begin{equation}
\label{cyclic_3}
E\Big((D_1^{\rho}-1)D_1^{\rho}s^{D_1^{\rho}}\Big)=
\left({1+\rho\over \rho}\right)s^2\left({\rho\over 1+\rho -s}\right)^{\rho+2}.
\end{equation}
It follows from (\ref{cyclic_1})-(\ref{cyclic_3}) and routine calculations that  for $1\leq k<n$,
$$
\Pr\{X_n^{\rho}=k\}={k\over n-k}\rho^k(1+\rho){\Gamma(n\rho)n!\over (n-k-1)!\Gamma(n\rho +k+1)}
$$
\begin{equation}
\label{cyclic_4}
=k\rho^k(1+\rho){(n)_{k}\over (n\rho +k)_{k+1}}
\end{equation}
and for $k=n$, we have
\begin{equation}
\label{cyclic_n}
\Pr\{X_n^{\rho}=n\}=\Pr\{\hat{D}^{\rho}_{i,n}=1, 1\leq i\leq n\}={\rho^n
n!\Gamma(n\rho)\over \Gamma(n+n\rho)}.
\end{equation}
\vskip .05in
Next , we consider the limiting distribution of $X_n^{\rho}$:
Fix $0<x<\infty$ and suppose that $k=\lfloor x\sqrt{n}\rfloor$, then we have
$$
\Pr\{X_n^{\rho}=k\}=k\rho^k(1+\rho){(n)_{k}\over (n\rho +k)_{k+1}}
$$
$$
={k\over n}\left({1+\rho\over \rho}\right){ (1-{1\over n})\cdot\cdot\cdot
(1-{k-1\over n})\over
(1+{k\over n\rho})\cdot\cdot\cdot (1+{1\over n\rho})}
$$
\begin{equation}
\label{cyclic_5}
\sim\left({1+\rho\over \rho}\right)x\exp\left(-{(1+\rho)x^2\over 2\rho}\right)
{1\over\sqrt{n}}.
\end{equation}
Hence $X_n^{\rho}/\sqrt{n}$ converges in distribution to a variable $\tilde{X}_{\rho}$ with
density
\begin{equation}
\label{cyclic_6}
f_{\tilde{X}_{\rho}}(x)=\left({1+\rho\over \rho}\right)x\exp\left(-{(1+\rho)x^2\over 2\rho}\right)
~\mbox{for }x\geq 0,
\end{equation}
and hence $E( X_n^{\rho})\sim
\sqrt{{\rho\pi n\over 2(1+\rho)}}$.
Let  $N_n^{\rho}$ denote the number of components in $G_n^{\rho}$,
then, using standard arguments (see \cite{Stepanov69}), it follows
from Corollary 3 and Theorem 4 that $(N_n^{\rho}-{1\over 2}\log n)/\sqrt{{1\over 2}
\log n}$ converges in distribution to the standard $N(0,1)$ distribution.
\par We apply  Theorem 3 and Theorem 4 to obtain the distribution
for $C_1^{\rho}(n)$:
$$
\Pr
\{C_1^{\rho}(n)=\ell\}=\Pr\{C_1^{\hat{D}(\rho,n)}=\ell\}
$$
\begin{equation}
\label{typical_1}
={\ell\over n}\Pr\{{\cal B}_{\ell}^{\hat{D}(\rho,\ell)}\}\Pr
\left\{\sum_{i=1}^{\ell} D_i^{\rho} \bigg| \sum_{i=1}^n D_i^{\rho}=n\right\}
\end{equation}
for $1\leq \ell\leq n$. To obtain a local limit theorem for the distribution
of $C_1^{\rho}(n)$, fix $0<x<\infty$ and suppose that $\ell=\lfloor xn \rfloor$. Then
it follows from Theorem 4, Corollary 2, (\ref{cyclic_4}), and (\ref{cyclic_n}) that
$$
\Pr\{{\cal B}_{\ell}^{\hat{D}(\rho,\ell)}\}=\sum_{k=1}^{\ell}{1\over k}
\Pr\{X_{\ell}^{\hat D(\rho,\ell)}=k\}
=\sum_{k=1}^{\ell-1}\rho^k(1+\rho){(\ell)_{k}\over (\ell\rho +k)_{k+1}}
+{\rho^n
\ell!\Gamma(\ell\rho)\over \Gamma(\ell+\ell\rho)}
$$
\begin{equation}
\label{typical_2}
\sim{1\over\sqrt{\ell}}\int_0^{\infty}
\left({1+\rho\over \rho}\right)\exp\left(-{(1+\rho)x^2\over 2\rho}\right) dx
 =\sqrt{{1+\rho\over \rho}}\cdot
\sqrt{{\pi\over 2\ell}}.
\end{equation}
Also, since the variables $D_1^{\rho}, D_2^{\rho}, ..$ are independent, we have
\begin{equation}
\label{typical_3}
\Pr
\left\{\sum_{i=1}^{\ell} D_i^{\rho} \bigg| \sum_{i=1}^n D_i^{\rho}=n\right\}
={\Pr\{\sum_{i=1}^{\ell}D_i^{\rho}\}\Pr\{\sum_{i=\ell+1}^n D_i^{\rho}\}
\over \Pr\{\sum_{i=1}^n D_i^{\rho}\}}
\end{equation}\begin{eqnarray*}
&=&{n\choose\ell}{\Gamma(\ell(1+\rho))\Gamma((n-\ell)(1+\rho))\Gamma(n\rho)\over
\Gamma(\ell\rho)\Gamma((n-\ell)\rho)\Gamma(n(1+\rho))}\\
&\sim& \sqrt{{\rho\over 1+\rho}}\sqrt{{n\over 2\pi\ell(n-\ell)}}.\\
\end{eqnarray*}
Substituting (\ref{typical_2}) and (\ref{typical_3}) into (\ref{typical_1}), we obtain
\begin{equation}
\label{typical_4}
\Pr
\{C_1^{\rho}(n)=\ell\}=
\Pr
\{C_1^{\rho}(n)=\lfloor xn\rfloor\}\sim {1\over 2n\sqrt{1-x}}.
\end{equation}
So,   ${C_1^{\rho}(n)\over n}$ converges in distribution to $Z_1$ as $n\to\infty$,  where $Z_1$ has density
$f_{Z_1}(u)=\theta(1-u)^{\theta-1}$ on the interval $(0,1)$ with $\theta=1/2$.
It follows from Theorem 3 and similar calculations that for
any integer $t\geq 1$ and constants $0<a_i<b_i<1$, where $1\leq i\leq t$,
$$
\lim_{n\to\infty}\Pr\left\{
a_i<{C_i^{\hat{D}(\rho,n)}\over n-C_1^{\hat{D}(\rho, n)}-\cdot\cdot\cdot
-C_{i-1}^{\hat{D}(\rho, n)}}<b_i, ~~1\leq i\leq t\right\}
$$
$$
=\prod_{i=1}^t\int_{a_i}^{b_i}{1\over 2\sqrt{1-x}}dx.
$$
 Hence, it follows from standard arguments (see, for example \cite{Hansen94}) that:
\begin{thm}
\label{Theorem4} The joint distribution of the normalized order statistics for the component sizes in $G_n^{\rho}$ converges to the {\it
Poisson-Dirichlet} (1/2) distribution on the simplex $\nabla=\{\{x_i\}:\sum x_i\leq 1, x_i\geq x_{i+1}\geq 0$ for every $i\geq 1\}$.
\end{thm}

\subsection{ An Anti-Preferential Attachment Model}
\vskip .1in
\noindent In this section we define $T_n^m:[n]\to[n]$, a random mapping with `anti-preferential
attachment', where $m\geq 1$ is a fixed interger parameter.
For $1\leq k\leq n$, we define $T_n^m(k)=Y_k^{(m,n)}$ where, as in the definition
of $T_n^{\rho}$, the variables $Y_1^{(m,n)}, Y_2^{(m,n)}, .., Y_n^{(m,n)}$ depend on the evolution of an urn scheme.
The distribution of each variable $Y_k^{(m,n)}$ is determined by a
(random) $n$-tuple of non-negative weights $\vec{b}(k)=(b_1(k), b_2(k), ..., b_n(k))$ where, for $1\leq j\leq n$, $b_j(k)$ is the number of balls
in the $j^{th}$ urn at the start of the $k^{th}$ round of the urn scheme.
Specifically, given $\vec{b}(k)=\vec{b}=(b_1, ..., b_n)$, we define
$$
\Pr\left\{Y_k^{(m,n)}=j \big| \vec{b}(k)=\vec{b}\right\}={b_j\over \sum_{i=1}^{n}b_i}.
$$
The random weight vectors  $\vec{b}(1), \vec{b}(2), ..., \vec{b}(n)$
associated with the urn scheme are determined
recursively. For $k=1$, we set $b_1(1)=b_2(1)=\cdot\cdot\cdot=b_n(1)=m$.
For $k>1$, $\vec{b}(k)$ depends on both $\vec{b}(k-1)$ and the value
of $Y_{k-1}^{(m,n)}$ as follows: Given that $Y_{k-1}^{(m,n)}=j$, we set $b_j(k)=b_j(k-1)-1$
and for all other $i\not=j$, we set $b_i(k)=b_i(k-1)$ (i.e.
if $Y_{k-1}^{(m,n)}=j$ then a ball is removed from the $j^{th}$ urn).
\par
The random mapping $T_n^m$ as defined above is an anti-preferential attachment model
in the following sense. Since, for $1\leq k\leq n$, we have $T_n^m(k)=Y_k^{(m,n)}$,
and
since the (conditional) distribution of $Y_k^{(m,n)}$ depends
on the state of the urn scheme at the start of round $k$,
it is clear that vertex $k$ is less likely to ` choose'  vertex $j$ if the weight $b_j(k)$
is (relatively) small, i.e. if several of the vertices $1,2,..., k-1$ have already been mapped to vertex $j$. It is also clear from the definition of $T_n^m$
that the in-degree of any vertex in the random digraph $G_n^m\equiv G(T_n^m)$
is at most $m$ and in the case $m=1$, $T_n^1$ is a (uniform) random permutation.
In the following proposition we determine the distribution of $T_n^m$.
\begin{thm}
\label{Theorem 6}
Suppose that $D^m_1, D^m_2, ...$ are i.i.d. $Bin(m, {1\over m})$
variables
where $m\geq 1$ is a fixed integer parameter.

\noindent
Let $\hat{D}(m,n)=({\hat D}_{1,n}^m, {\hat D}_{2,n}^m, ...,
{\hat D}_{n,n}^m)$ be a sequence of variables  with joint distribution given by
$$
\Pr\{{\hat D}_{i,n}^m=d_i, 1\leq i\leq n\}=\Pr\left\{ D_i^m=d_i, 1\leq i\leq n\bigg|
\sum_{i=1}^n D_i^m= n\right\}.
$$
 Then the random mappings $T_n^m$ and $T_n^{\hat D(m,n)}$ have the same distribution.

\end{thm}
\begin{proof}
To prove the result it is enough to show that for any $n\geq 1$ and any  $f\in{\cal M}_n$
$$
\Pr\{T_n^m=f\}=\Pr\{T_n^{\hat D(m,n)} =f\}.
$$
Suppose that $f\in{\cal M}_n$ and that $\vec{d}(f)=(d_1, d_2, ..., d_n)$.
It is straightforward to check that
$$
\Pr\{{\hat D}_{i,n}^m=d_i, 1\leq i\leq n\}=\Pr\left\{ D_i^m=d_i, 1\leq i\leq n\bigg|
\sum_{i=1}^n D_i^m= n\right\}
$$
$$
={\prod_{i=1}^n{m\choose d_i}\over
{nm\choose n}},
$$
and hence, from the definition of $T_n^{\hat D(m,n)}$, that
$$
\Pr\{T_n^{\hat D(m,n)}=f\} = {d_1!\cdot\cdot\cdot d_n!\over n!}\times
{\prod_{i=1}^n{m\choose d_i}\over
{nm\choose n}}
$$
$$
={\prod_{i=1}^n(m)_{d_i}\over (nm)_n}.
$$
On the other hand,  $T_n^m=f$  with $\vec{d}(f)=(d_1, d_2, ..., d_n)$ if and only if
for each $1\leq i\leq n$, $d_i$ balls are removed, in a certain order,
from the $i^{th}$ urn during the evolution of the urn model
described above.
So it follows from the definition of $T_n^m$ in terms of the
urn scheme,
that
$$
\Pr\{T_n^m=f\}={\prod_{i=1}^n(m)_{d_i}\over (nm)_n}.
$$
The result follows since $n\geq 1$ and $f\in{\cal M}_n$ were arbitrary.
\end{proof}

Again, it is clear from Theorem 6 that the order in which a realisation of
$T_n^{m}$ is sequentially constructed does not matter.
In particular, suppose that $i_1, i_2, ..., i_n$ is a permutation of $[n]$
and for $1\leq k\leq n$, let $\tilde{T}_n^{m}(i_k)\equiv Y_n^{(m, n)}$
where the variables $Y_1^{(m, n)}, Y_2^{(m, n)}, ..., Y_n^{(m, n)}$
are as defined above. Then it follows from the proof of Theorem 6
that $\tilde{T}_n^{m}\stackrel{d}{ \sim}T_n^{\hat{D}(m,n)}\stackrel{d}{ \sim} T_n^{m}$.
\par We apply Theorem 6 to investigate the distributions
of the number of cyclic vertices in $G_n^m\equiv G(T_n^m)$
and the size of a typical component in $G(T_n^m)$.
Let $X_n^m$
denote the number of cyclic vertices
in the random digraph $G_n^m\equiv G(T_n^m)$ and let
$C_1^m(n)$ denote the size of the component in $G_n^m$
which contains the vertex 1.
Since $T_n^m\stackrel{d}{ \sim} T_n^{\hat D(m,n)}$, we have $G(T_n^m)\stackrel{d}{ \sim} G_n^{\hat D(m,n)}$
and $X_n^m\stackrel{d}{ \sim} X_{n}^{\hat D(m,n)}$. So, it follows from Theorem 1 and
Theorem 6 (and its proof) that for $m\geq 2$ and $1\leq k<n$, we have
$$
\Pr\{X_n^m=k\}=\Pr\{X_n^{\hat D(m,n)}=k\}=
{k\over n-k}E((\hat D^m_{1,n}-1)\hat D^m_{1,n}\hat D^m_{2,n}\cdot
\cdot\cdot \hat D_{k,n})
$$
\begin{eqnarray*}
&=&{k\over n-k}\sum_{\vec d~ s.t. \sum_{i=1}^nd_i=n}
(d_1-1)d_1d_2\cdot\cdot\cdot d_k\times{{m\choose d_1}\cdot\cdot\cdot
{m\choose d_n}\over {nm\choose n}}\\
&=&{k\over n-k}\sum_{t=k+1}^{\min(n,km)}\sum_{\vec d~ s.t. \sum_{i=1}^kd_i=t\atop
and ~\sum_{i=1}^nd_i=n}
(d_1-1)d_1d_2\cdot\cdot\cdot d_k\times{{m\choose d_1}\cdot\cdot\cdot
{m\choose d_n}\over {nm\choose n}}\\
\end{eqnarray*}
\begin{eqnarray*}
&=&{k\over n-k}\sum_{t=k+1}^{\min(n,km)}\sum_{\vec d~ s.t. \sum_{i=1}^kd_i=t}
(d_1-1)d_1d_2\cdot\cdot\cdot d_k\times{{m\choose d_1}\cdot\cdot\cdot
{m\choose d_k}{nm-km\choose n-t}\over {nm\choose n}}\\
&=&{k\over n-k}m^k(m-1)\sum_{t=k+1}^{\min(n,km)}\sum_{\vec d~ s.t. \sum_{i=1}^kd_i=t}
{{m-2\choose d_1-2}{m-1\choose d_2-1}\cdot\cdot\cdot
{m-1\choose d_k-1}{nm-km\choose n-t}\over {nm\choose n}}\\
&=&{k\over n-k}m^k(m-1)\sum_{t=k+1}^{\min(n,km)}
{{km-k-1\choose t-k-1}{nm-km\choose n-t}\over {nm\choose n}}\\
&=&{k\over n-k}m^k(m-1)
{{nm-k-1\choose n-k-1}\over {nm\choose n}}.\\
\end{eqnarray*}
 In the summations above the sum is always taken over those degree sequences
for which the binomial coefficients are defined. We also adopt the formal convention that
${0\choose 0}=1$.
 Finally, for $k=n$ and $m\geq 2$, we obtain
$$\Pr\{X_n^m=n\}=\Pr\{\hat D_{i,n}^m=1,1\leq i\leq n\}={m^n\over {nm\choose n}},
$$
and
when $m=1$, we have  $X_n^1\equiv n$.

To obtain a local limit theorem for $X_n^m$,  fix $0<x<\infty$ and suppose that $k=\lfloor x\sqrt{n}\rfloor$.
Then we have
$$
\Pr\{X_n^m=k\}
={k\over n-k}m^k(m-1)
{{nm-k-1\choose n-k-1}\over {nm\choose n}}
$$
$$
={k\over n-k}m^k(m-1)
{(n)_{k+1}\over (nm)_{k+1}}
$$
$$
\sim \left({m-1\over m}\right)x\exp\left({-(m-1)x^2\over 2m}\right){1\over\sqrt{n}}.
$$
It follows that $X_n^m/\sqrt{n}$ converges in distribution to a variable $\tilde{X}_m$ with  density
$$
f_{\tilde{X}_m}(x)=\left({m-1\over m}\right)x\exp\left({-(m-1)x^2\over 2m}\right)\quad\mbox{for }
x\geq 0,
$$
and hence $E(X_n^{m})\sim\sqrt{{mn\pi\over 2(m-1)}}$. Let $N_n^m$ denote the number of components
in $G_n^m$. Then, by  the same argument  as used for $N_n^{\rho}$ above, it follows that\\
$(N_n^m-{1\over 2}\log n)/\sqrt{{1\over 2}\log n}$ converges in distribution to a standard $N(0,1)$
distribution.
\par We apply  Theorem 3 and Theorem 6 to obtain the distribution
for $C_1^{m}(n)$:
$$
\Pr
\{C_1^{m}(n)=\ell\}=\Pr\{C_1^{\hat{D}(m,n)}=\ell\}
$$
\begin{equation}
\label{typical_anti1}
={\ell\over n}\Pr\{{\cal B}_{\ell}^{\hat{D}(m,\ell)}\}\Pr
\left\{\sum_{i=1}^{\ell} D_i^{m} \bigg| \sum_{i=1}^n D_i^{m}=n\right\}
\end{equation}
for $1\leq \ell\leq n$. To obtain a local limit theorem for the distribution
of $C_1^{m}(n)$, fix $0<x<\infty$ and suppose that $\ell=\lfloor xn \rfloor$. Then
it follows from Theorem 6  and Corollary 2,  that
$$
\Pr\{{\cal B}_{\ell}^{\hat{D}(m,\ell)}\}=\sum_{k=1}^{\ell}{1\over k}
\Pr\{X_{\ell}^{\hat D(m,\ell)}=k\}
$$
$$
=\sum_{k=1}^{\ell-1}{m^k(m-1)\over \ell-k}{{\ell m-k-1\choose \ell-k-1}\over
{\ell m\choose \ell}}\quad+\quad{1\over \ell}{m^{\ell}\over {\ell m\choose \ell}}
$$
\begin{equation}
\label{typical_anti2a}
\sim{1\over\sqrt{\ell}}\int_0^{\infty}\left({m-1\over m}\right)\exp
\left(-{(m-1)x^2\over 2m}\right) dx =\sqrt{{m-1\over m}}
\sqrt{{\pi\over 2\ell}}   \, .
\end{equation}
Since $D_1^m, D_2^m, ... $ are i.i.d. $Bin(m, {1\over m})$ variables, we also
have
\begin{equation}
\label{typical_anti3}
\Pr\left\{\sum_{i=1}^{\ell} D_i^m=\ell
\Big|\sum_{i=1}^n D_i^m=n\right\}
={{m\ell\choose\ell}{mn-m\ell\choose
n-\ell}\over {mn\choose n}}
\sim {\sqrt{mn}\over \sqrt{2\pi\ell(m-1)(n-\ell)}}.
\end{equation}
Substituting (\ref{typical_anti2a}) and (\ref{typical_anti3}) into (\ref{typical_anti1})
we obtain
\begin{equation}
\label{typical_anti2}
\Pr\{C_1^m(n)=\ell\}
\sim{1\over 2n\sqrt{1-x}}\quad\mbox{as }n\to\infty.
\end{equation}
 It follows that as $n\to\infty$, ${C_1^m(n)\over n}$ converges in distribution to $Z_1$, where $Z_1$ has
$f_{Z_1}(u)=\theta(1-u)^{\theta-1}$ on the interval $(0,1)$ with parameter $\theta=1/2$. Again, as in the case of the preferential attachment
model, it is not difficult to extend (\ref{typical_anti2}) in order to show that:

\begin{thm}
\label{Theorem5} The joint distribution of the normalized order statistics for the component sizes in $G_n^m$ converges to the {\it
Poisson-Dirichlet} (1/2) distribution on the simplex $\nabla=\{\{x_i\}:\sum x_i\leq 1, x_i\geq x_{i+1}\geq 0$ for every $i\geq 1\}$.
\end{thm}%
%

\section{Final Remarks }
\label {s:final}

In this paper we have introduced a new random mapping model, $T_n^{\hat{D}}$,
which is defined in terms of a collection of exchangeable, non-negative, integer-valued
`in-degree'  variables $\hat{D}_1, \hat{D}_2, ..., \hat{D}_n$ such that $\sum_{i=1}^n
\hat{D}_i=n$. We note the classic model $T_{{\bf p}(n)}$ has the property that each vertex
{\it independently} chooses its image under $T_{{\bf p}(n)}$. In
$T_n^{\hat D}$ we have replaced this independence with
 exchangeability  of the in-degrees in $G_n^{\hat D}$ 
(and in the examples with the
`independence' of random variables which are used to define the in-degrees).

We have shown that the joint distribution of the variables
$\hat{D}_1, \hat{D}_2,...,\hat{D}_n$ is the key to understanding both the local and
global structure of  $G_n^{\hat D}$, and
 that the distributions of several variables associated with the structure of $G_n^{\hat D}$
 can be
expressed in terms of  expectations of various functions of $\hat{D}_1, \hat{D}_2, ..., \hat{D}_n$.
In the special case where the variables $\hat{D}_1, \hat{D}_2, ..., \hat{D}_n$ have the same
distribution as a collection of i.i.d. variables $D_1, D_2, ..., D_n$
{\it conditioned} on $\sum_{i=1}^n D_i=n$,  we have shown that it is easy
to apply our results to obtain exact and asymptotic distributions for the number of cyclic
vertices, the number of components, and the size of a typical component in $G_n^{\hat D}$.
Both $T_{{\bf p}(n)}$ and $T_n^{\hat D}$ share the property that these mappings restricted
to the cyclical vertices are both uniform random permutations on the cyclical vertices.
This means that given the distribution of the number of cyclic vertices, $X_n^{\hat D}$,
 (see Theorem 1) we can exploit known results for uniform random permutations
 to obtain results for $T_n^{\hat D}$, and it explains
 some similarities between
 results for both $T_{{\bf p}(n)}$ and
 $T_n^{\hat D}$\,.

The transparent relationship between the distribution of the variables $\hat D_1, \hat D_2, ...,\hat D_n$
and the structure of random mapping $T_n^{\hat D}$ is the main advantage of this new
model. It also provides the basis for testing  and fitting random mapping models in various
applications. For example, if $\hat{D}_1, \hat{D}_2, ..., \hat{D}_n$ have the same
distribution as $n$ independent {\it Poisson}(1) variables, $D_1, D_2, ..$
$.., D_n$,
{\it conditioned} on $\sum_{i=1}^n D_i=n$, then $T_n^{\hat D}$ is a uniform
random mapping and, for large $n$ and $1\leq i\leq n$, the marginal distribution
of each $\hat{D}_i$ is approximately {\it Poisson}(1).  So vertex in-degree
data, such as directed epidemic contacts, can be used to test whether the
associated mapping is a realisation of a uniform random mapping model.

As another example, we mention the work of Arney and Bender on random
mappings with constraints on coalescence \cite{Arney82}. Their work was motivated, in part, by the analysis
of shift register data. In order to model a random shift register
they put a uniform measure on ${\cal M}_n^{\{0,1,2\}}$, the set of all mappings
$f:[n]\to [n]$ such that, for every $1\leq i\leq n$, $|f^{-1}(i)|$, the number of
pre-images of $i$ under $f$, equals $0,1,$ or $2$. So, if $f\in{\cal M}_n^{\{0,1,2\}}$,
then every vertex in $G_n(f)$ has in-degree equal $0,1,$ or $2$. Arney and Bender
observed that in some respects their model does not fit the shift register data.
In particular, their model predicts $0.293 n$ vertices with in-degree $0$ whereas
the average number of vertices with in-degree $0$ in a random shift register
is $n/4$. By using the model $T_n^{\hat D}$ instead, we can more successfully capture
the local structure of the shift register data. Specifically, suppose that $\hat{D}_1, \hat{D}_2,
..., \hat{D}_n$ have the same distribution as $n$ independent $Bin(2, {1\over 2})$ variables,
$D_1, D_2, ..., D_n$, conditioned on $\sum_{i=1}^n D_i=n$.
Then $\Pr\{\hat{D}_1 =0\}={1\over 4}(1+{1\over 2n-2})^{-1}$ and the expected number of
vertices with in-degree 0 in $T_n^{\hat D}$ is asymptotic to ${n\over 4}$.
We note that, for example, the asymptotic distribution of the
normalised typical
component size is the {\it same} for both the Arney and Bender model and $T_n^{\hat D}$.
So it is not surprising that Arney and Bender found that their model fit some
other features of the shift register data quite well.

Finally, we mention that in applications of random mapping models  in cryptology
and in epidemic process modelling  the distributions of the number of predecessors
and of the number of successors of an arbitrary vertex or set of vertices are also of interest.
We develop a calculus for computing these distributions based on the underlying
variables $\hat{D}_1, \hat{D}_2, ..., \hat{D}_n$ in a companion paper \cite{HJ06a}.


\end{document}